\documentclass[final]{article}
\usepackage{showkeys}
\usepackage{amsmath, amssymb, epsf, a4wide,amssymb,amsmath,amsthm}
\usepackage{euscript}
\usepackage{graphics,graphicx}

\newcommand{\N}{\Bbb N}
\newcommand{\R}{\Bbb R}

\newcommand{\Q}{\Bbb Q}

\mathchardef\flat="115B

\newtheorem{thm}{Theorem}[section]
\newtheorem{lemma}[thm]{Lemma}

\newtheorem{prop}[thm]{Proposition}

\newtheorem{prob}[thm]{Problem}

\begin{document}

\vskip 2.7cm

\centerline{\large{\bf MANIFOLDS: HAUSDORFFNESS VERSUS
HOMOGENEITY}}

\vskip 0.5cm

\centerline{ \normalsize MATHIEU BAILLIF AND ALEXANDRE GABARD}

\vskip 1cm

\smallskip
\newbox\abstract
\setbox\abstract\vtop{\hsize 11 cm \noindent

\footnotesize \noindent\textsc{Abstract.} We analyze the
relationship between Hausdorffness and homogeneity in the frame of
manifolds, not confined to be Hausdorff. We exhibit examples of
homogeneous non-Hausdorff manifolds and prove that a Lindel\"of
homogeneous manifold is Hausdorff. }

\centerline{\hbox{\copy\abstract}}

\bigskip



\normalsize

\hskip 1.5cm

\section{Introduction}\label{sec1}

Our purpose here is to analyze the relationship between
Hausdorffness and homogeneity in the frame of manifolds. We
give the word {\it manifold} its broadest sense, that is, a
topological space locally homeomorphic to the Euclidean space
${\Bbb R}^n$ of a fixed dimension (without assuming the Hausdorff
separation axiom).

Recall that a connected Hausdorff manifold $M$ is {\it
homogeneous}, i.e. for each $x,y \in M$, there is a homeomorphism
$h: M \to M$ taking $x$ to $y$ (see \cite{vanDantzig} or
\cite{Vick}, p.$\,$150).

This property is true only under the Hausdorff assumption. Without
it, one may well have an non homogeneous manifold, for example the
well known {\it line with two origins}: take two copies of the
real line ${\Bbb R}$ and identify all corresponding points of the
copies but the origin (Figure 1). This yields a one-dimensional
manifold in which the two origins cannot be separated\footnote{We
say that two points of a topological space can be {\it separated}
if there are two disjoint open sets containing one of them each.}.
Notice though that a point different from the origins can be
separated from any other point, so the manifold is not
homogeneous. Another well known example of non-Hausdorff manifold
is the {\it branching line} obtained by identifying the points
$<0$ in the two copies of $\R$ (Figure 1).

\begin{figure}[h]
\centering
    \includegraphics[angle=0]{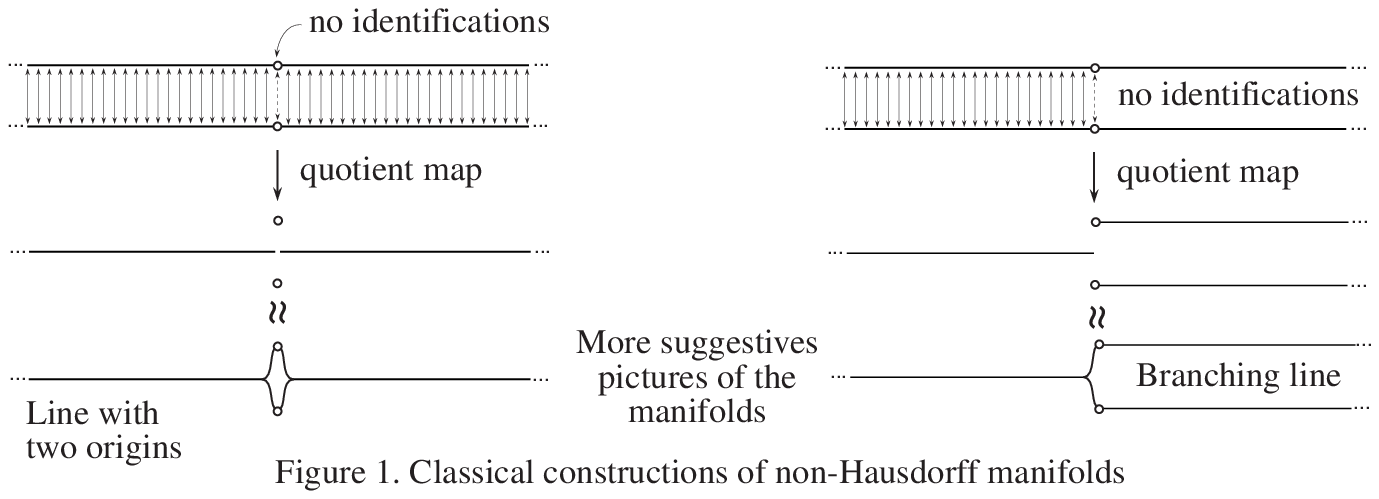}
\end{figure}

One may think that homogeneity is a sufficient condition to
characterize, in the realm of manifolds, those which are
Hausdorff. We show that this is not the case by exhibiting two
examples. The first, called the {\it complete feather} or {\it
everywhere branching line} $F$ will be discussed in $\S$
\ref{sec2}. It was first defined by Haefliger and Reeb in
\cite{HaefligerReeb}, and is constructed by ``grafting'' lines to
all points of a line, and iterating this process indefinitely (see
Figure 2).

\begin{figure}[h]
\centering
    \includegraphics[angle=0]{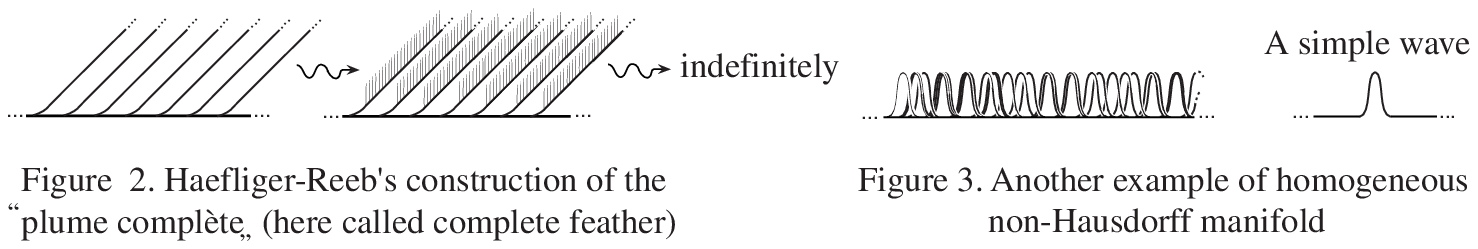}
\end{figure}

$F$ is a non-Hausdorff homogeneous $1$-manifold but is neither
separable\footnote{A space is {\it separable} if it has a
countable dense subspace.} nor Lindel\"of\footnote{A space is {\it
Lindel\"of} if each open cover has a countable subcover.}. It
furthermore has some interesting contractibility properties: it is
contractible but not strongly (i.e. in such a way that the
collapsing point stays fixed during the deformation), even though
each point has such a strongly contractible neighborhood. This
answers a question left open in \cite{Gauld:1983}.

Our second example is the {\it everywhere doubled line} $D$, a
``continuous'' version of the line with two origins in which we
perform the duplication process at all points (Figure 3 tries to
give a representation). $D$ is homogeneous, separable but neither
Hausdorff nor Lindel\"of. A discussion of this example is the
object of $\S$ \ref{sec3}. The quest for a Lindel\"of example is
ruled out by the following:

\begin{thm}\label{thm1}
A homogeneous Lindel\"of manifold is Hausdorff (and therefore
metrizable\footnote{Recall that a Lindel\"of Hausdorff manifold is
metrizable; this follows from Urysohn's metrization theorem, since
Lindel\"of and locally second countable imply second countable.
}).
\end{thm}

This will be proved in $\S$ \ref{sec4}. The proof uses that a
manifold is a Baire space\footnote{Every countable intersection of
dense open sets is dense.} (since a space that is locally
Baire\footnote{For
  any property $\cal P$ attributable to a space, a
  space is said to be {\it locally} ${\cal P}$ if each of its
  points has a neighborhood with the property $\cal P$.} is in fact Baire). Alternatively
  we can also argue that a (non necessarily Hausdorff) locally compact
space is Baire, a fact that is usually proved only for Hausdorff
spaces, see for instance \cite{Bourbaki}, but which remains valid
in this more general setting (see $\S$ \ref{sec5}).

\section{The complete feather $F$}\label{sec2}

Let us first give a loose description of $F$. The idea is to start
with the usual real line, and to add branches (like in the
branching line) at any $x\in\R$. This results in a ``hairy line'',
with branches at level $1$. Then, we continue the process by
adding new branches at level $2$ to all points in branches in
level $1$, and so on indefinitely. The resulting space $F$ is a
$1$-manifold whose homogeneity comes from the fact that we can
``flip'' a branch at level $i$ with a branch at level $i+1$. It is
also contractible.

Now, the formal definition.
The underlying set of $F$ is
$$
  \left\{ s\,:\,
    \begin{array}{l} s=(s_0,\dots,s_n)\text{ for some }\,n\ge 0, \,s_i\in\R,\\
    \text{ and }
    s_0<s_1<\dots < s_{n-1} \le s_n \end{array}
  \right\} .
$$
Notice that the last inequality is {\em not} strict. One should
interpret the sequences of length $1$ as the usual real line,
those of length $\le 2$ as the hairy line, and so on. We
topologize $F$ with the order topology for the following partial
order:
$$
(s_0, \dots,s_n)<(t_0,\dots,t_m) \text{  iff  } n\le m, \;\;
s_i=t_i \text{ for } i=0,\dots,n-1 \text{ and } s_n<t_n .
$$

Notice that $(s_0, \dots,s_n)$ and $(s_0, \dots,s_n,s_n)$ are
incomparable and have the same predecessors. $F$ is a ``tree'' in
the sense that the predecessors of any point are totally ordered.

\begin{figure}[h]
\centering
    \includegraphics[angle=0]{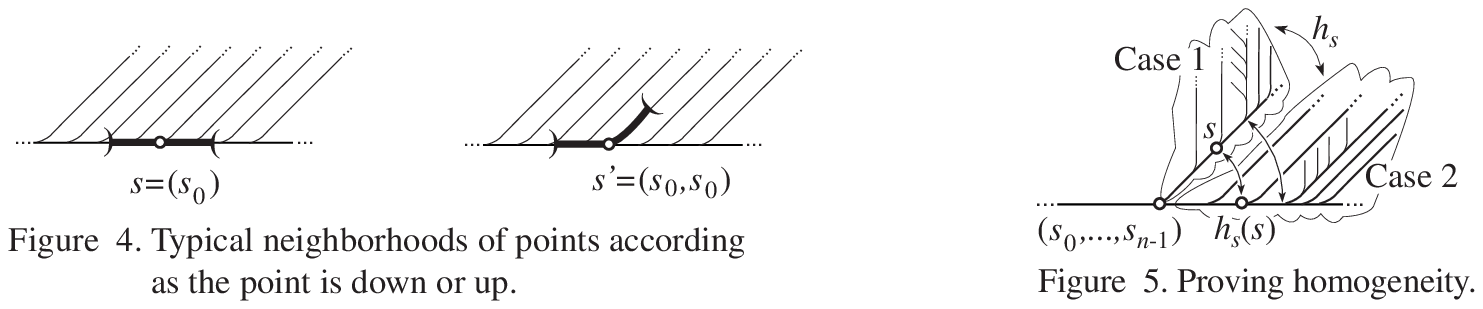}
\end{figure}

\begin{prop}
  $F$ is a connected homogeneous non-Hausdorff $1$-manifold.
\end{prop}
Notice that $F$ is non separable, in fact there is even an
uncountable family of pairwise disjoint open sets in $F$. \proof
  One sees immediately that $F$ is a non-Hausdorff manifold, since the intervals are homeomorphic to $\R$ and
  the points of the form $(s_0,\dots,s_n)$ and $(s_0,\dots,s_n,s_n)$ cannot be separated.
  Connectedness is also easy.
  To see that $F$ is homogeneous, we first show that given any point $s=(s_0,\dots,s_n)\in F$,
  there is an homeomorphism of $F$ sending $s$ to $(s_n)\in F$.
  We first consider the map $h_s:F \to F$ that flips the two ``branches''
  emanating from $(s_0,\dots,s_{n-1})$ (see Figure 5). It is given
  by the formula
  $$
    h_s(r)\!=\!\!
    \left\{
    \begin{array}{cl}
    (s_0,...,s_{n-2},r_n,r_{n+1},...,r_m)&\!\!\!\!\text{ if } r=(s_0,...,s_{n-2},s_{n-1},r_n,r_{n+1},...,r_m), \\
    \!\!\! (s_0,...,s_{n-2},s_{n-1},r_{n-1},r_n,...,r_m)&\!\!\!\!\text{ if } r=(s_0,...,s_{n-2},r_{n-1},r_n,...,r_m)
     \text{ with } r_{n-1}\ge s_{n-1}, \\
    r &\!\!\!\!\text{ otherwise}.
    \end{array}
    \right.
  $$
One sees easily that $h_s$ is an homeomorphism (actually it is an
involution).
   For $s=(s_0,\dots,s_n)\in F$ and $k\le n$, let $s^{(k)}=(s_0,\dots,s_{n-k})$.
   Then, $h_{s^{(n-1)}}\circ\cdots\circ h_{s^{(1)}}\circ h_s(s)=(s_n)$.
   To finish the proof it suffices to remark that for $t\in\R$ the map
   $(s_0,\dots,s_n)\mapsto(s_0+t,\dots,s_n+t)$ is a homeomorphism.
\endproof

\begin{lemma}
  $F$ is contractible.
\end{lemma}
\proof
  The idea is to contract all points of the form $(s_0,\dots,s_{n-1},s_n)$ on the point
  $(s_0,\dots,s_{n-1},s_{n-1})$ between time $\frac{1}{n+1}$ and $\frac{1}{n}$.
  For $x\in\R$ and $t\in [0,1]$, let $\varphi_t^x:{\Bbb R} \to {\Bbb
  R}$ be defined by
  $$
    \varphi_t^x(y)=\left\{
    \begin{array}{cl} y&\text{ if } y<x\\
    (1-t)y+tx&\text{ if } y\ge x .
    \end{array}
    \right.
  $$
  Then, $\varphi_0^x=\text{\rm id}$, and $\varphi_1^x(y)=x$ for all $y\ge x$.
  If $s=(s_0,\dots,s_n)\in F$, we let
  $\phi_t(s)= (s_0,\dots,s_{n-1},\varphi_t^{s_{n-1}}(s_n))$, so
  $\phi_1(s)=(s_0,\dots,s_{n-1},s_{n-1})$.
  If $s=(s_0,\dots,s_n)$, we write $t'=n(n+1)t-n$ and define
  $$
    h_t(s)=\left\{
    \begin{array}{cl}
      s&\text{ if } n=0 \text{ or } t\le\frac{1}{n+1},\\
      \phi_{t'}(s) &\text{ if } n\neq 0 \text{ and } t\in [\frac{1}{n+1},\frac{1}{n}],\\
      h_t((s_0,\dots,s_{n-1})) & \text{ if } n\neq 0 \text{ and } t>\frac{1}{n}.
    \end{array}
    \right.
  $$
The definition is implicit, but this causes no problem: we proceed
by induction on $n$. Thus, $h_0=\text{\rm id}$ and
$h_1(s)=(s_0,s_0)$ if $s=(s_0,\dots,s_n)$ with $n>0$,
$h_1((s_0))=(s_0)$. We then define $h_t((s_0))=h_t((s_0,s_0))=
(s_0-t+1)$ for $t\in[1,2]$. It is not difficult to see that $h_t$
is continuous and that $h_2(F)$ is included in the sequences of
length $1$ which are homeomorphic to $\R$ and thus contractible.
\endproof

A space $X$ is {\it strongly contractible to the point $p$} if
there exists an homotopy $h_{t}:X \to X$ such that
$h_0=\text{id}$, $h_1\equiv p$ and $h_t(p)=p$ for all $t$. D.
Gauld \cite{Gauld:1983} showed that if $X$ is contractible,
locally strongly contractible to $p$ and completely regular at
$p$, then $X$ is strongly contractible to $p$. Further he asked
whether ``completely regular'' could be dropped; the complete
feather $F$ gives a counterexample since:

\begin{lemma}
  {$F$ is not strongly contractible to any of its points.}
\end{lemma}

\proof
  Call {\it twins} the pairs of points of $F$ of the form
  $\{(s_0,\dots,s_{n}), (s_0,\dots,s_{n},s_{n})\}$.
  Any sequence $(s_0,\dots,s_{n-1},s_{n}^m)$ $(m\in {\Bbb N})$ with $s_{n}^m\nearrow s_{n}$ converges
  to both twins. Thus, if one of the twins moves, the other must also move.
  Since any point of $F$ as a twin the result follows.
\endproof

\normalsize

\section{The everywhere doubled line $D$}\label{sec3}
\medskip

We can build $D$ either as an inductive limit or with two copies
of the line with an exotic topology. We give the latter
construction. The underlying set of $D$ is $\R \times \{0,1 \}$.
Points of $D$ with zero second coordinate are said ``down'', the
others ``up''. A base for the topology is given by usual open sets
downstairs with a finite (eventually zero) number of points
removed and lifted upstairs; that is, subsets of the form
$$
  \mathcal{U}_{O,F}=(O\backslash F)\times \{0 \} \cup F \times \{ 1 \},
$$
where $O\subset\R$ is open and $F\subset O$ finite. Such subsets
will be called {\it waves}. It is immediate that the waves are
closed under finite intersections, we topologize $D$ with the
topology given by this base.

\medskip
\begin{prop} $D$ is a connected
non-Hausdorff homogeneous separable 1-manifold.
\end{prop}

\proof First, it is clear that $D$ is non-Hausdorff, since two
points $x,y$ having the same first coordinate cannot be separated.
It is also immediate that $\Q\times\{0\}$ is dense in $D$, and
that any wave $\mathcal{U}_{O,F}$ where $O=\,]a,b[$ is
homeomorphic to $\R$. This proves that $X$ is a separable
1-manifold. Connectedness is also easy.

To check homogeneity, we begin by observing that the translations
$t_s: (t,i) \mapsto (t+s,i)$ with $s\in {\Bbb R}$ are clearly
homeomorphisms. This settles the case where the two given points
have the same second coordinates. If not, use the map exchanging
up and down at one value of the abscissaes
$$
  e_s: (t,i) \mapsto
    \left\{
    \begin{array}{ll}
       (t,1-i)&\text{ if }t=s \\
       (t,i)  &\text{ if }t\not= s,
    \end{array}
    \right.
$$
which is a homeomorphism, since it acts simply by adding or
removing an oscillation to a given wave (or eventually do nothing
at all if $s$ is outside from wave's range). \hfill $\square$

\medskip \small
{\it Note.} Working with reflections and exchange maps, we even see that $X$ is {\it involutorially
homogeneous} (i.e. the homeomorphism taking $x$ to $y$ can always
be chosen to be an involution).

\medskip \normalsize
{\it Remarks.} As a variant of this construction we can also
triple each points of the line. This gives a counterexample to an
erroneous claim made by Fuks-Rokhlin who asserted, that any
one-dimensional manifold becomes disconnected after ones removes
two suitably chosen points (see \cite{FuksRokhlin}, p.$\,$135).

It is a theorem of Nyikos \cite{Nyikos:1984} that a Hausdorff
connected manifold of dimension $\ge 1$ has cardinality the
continuum, the same construction starting with ${\Bbb R}\times
\kappa$ for $\kappa$ any cardinal, shows that (homogeneous)
connected manifolds can have arbitrarily large cardinality.

\section{Proof of Theorem \ref{thm1}}\label{sec4}

We will prove the following:

\begin{thm}\label{thm2}
  Let $X$ be a topological space which is
  homogeneous, Lindel\"of, locally Hausdorff and
  Baire. Then $X$ is Hausdorff.
\end{thm}
\medskip

Theorem \ref{thm1} is then immediate.
We will need the following application of
Zorn's Lemma:

\medskip
\begin{lemma}\label{lemmazorn}Let $X$ be a locally Hausdorff space. Then
  for each point $x\in X$, there exists a Hausdorff dense open
  set $U_x$ containing $x$.
\end{lemma}

\noindent{\it Proof of \ref{lemmazorn}.} Let $x \in X$. Consider
$\frak O_{x}$ the set of all Hausdorff open sets containing $x$,
ordered by inclusion. Since $X$ is locally Hausdorff, $\frak
O_{x}$ is non-empty. We check that $\frak O_{x}$ is inductive. Let
$\frak C$ be a totally ordered subset of $\frak O_{x}$. As usual
let $V=\cup_{U\in \frak C} U$ be the natural upper bound of $\frak
C$. Then $V$ is open, and Hausdorff: given two points in $V$ (say
$y,z$), each of them belongs to some $V_y,V_z \in \frak C$. But
since $\frak C$ is totally ordered both points belong to one of
them (say $V_y$), and can thus be separated by open sets of $V_y$,
which are also open in $X$.

By Zorn's lemma, there is $U_x$ maximal in $\frak O_{x}$. We check
its density. So, let $\mit\Omega$ be a non-empty open set of $X$.
We can assume that $\mit\Omega$ is Hausdorff. If $\mit\Omega\cap
U_x=\emptyset$, since $\mit\Omega$ and $U_x$ are Hausdorff, so is
their union, and thus $U_x\cup \mit\Omega\in\frak O_{x}$, which
contradicts the maximality of $U_x$.  \hfill $\square$
\medskip

\noindent{\it Proof of \ref{thm2}.} By homogeneity, it is enough
to prove the existence of a point $x_0\in X$ which can be
separated from each other point $y\in X$.

For all $x\in X$ let $U_x\ni x$ be given by Lemma \ref{lemmazorn}.
The collection $(U_x)_{x\in X}$ is an open cover of $X$ from which
we extract a countable subcover $(U_{x_i})_{i\in {\Bbb N}}$ (by
Lindel\"ofness). Since the $U_{x_i}$ are dense open sets and $X$
is Baire, their intersection $\cap_{i\in {\Bbb N}} U_{x_i}$ is
dense, and so in particular non-empty. Any point $x_0$ in this
intersection is separable from any other $y\in X$: since the
$(U_{x_i})_{i\in {\Bbb N}}$ cover $X$, $y$ is in $U_{x_i}$ for
some $i\in\N$; but so does $x_0$, and since $U_{x_i}$ is
Hausdorff, $x_0$ and $y$ can be separated. \hfill $\square$

\medskip
{\it Remarks.} The preceding results raise some problems we found
worth mentioning here. Firstly, the Lindel\"of condition in
Theorem \ref{thm1} is in a sense too strong, since it implies
metrizability, and there are non-metrizable Hausdorff manifolds.
\begin{prob}\label{prob1}
  In Theorem \ref{thm1}, can Lindel\"of be replaced by a weaker condition in order to
  ensure the Hausdorffness but not necessarily the metrizability of the manifold?
\end{prob}
Secondly, the homogeneous non-Hausdorff manifolds $F$ and $D$ both
contain an uncountable (closed) discrete subset: Take one point in
each branch at level one in $F$, and all the ``up'' points in $D$.
So, another problem is:
\begin{prob}\label{prob2}
   Is there a homogeneous non-Hausdorff manifold that contains no uncountable (closed) discrete
   subset, or even stronger that is hereditarily separable\footnote{In short $HS$, and means that every subspace is
separable.}?
\end{prob}
(Whether there are non-metrizable $HS$ {\em Hausdorff} manifolds
or not is known to be independent of ZFC. Under $CH$, Rudin-Zenor
\cite{RudinZenor} constructed a non-metrizable $HS$ manifold. On
the other hand Szentmikl\'ossy \cite{Szentmiklossy} showed that
under $MA+\neg CH$, every locally compact $HS$ Hausdorff space is
$HL$(=hereditarily Lindel\"of), and so metrizable if a manifold.
The result in \cite{Szentmiklossy} is actually stated with compact
instead of locally compact, but the above follows by taking the
one-point compactification.) A negative answer to the first part
of Problem \ref{prob2} would yield that
$\omega_1$-compact\footnote{A discrete closed set is at most
countable.} is an answer to Problem \ref{prob1}.

\section{About the Baire property}\label{sec5}

We call a space {\it quasi-compact} if from any open cover one may
extract a finite subcover, and {\it compact} if moreover it is
Hausdorff.
\begin{thm}[Baire slightly extended]\label{thm3}
   Let $X$ be a locally
   compact (not necessarily Hausdorff) space. Then $X$ is a Baire
   space.
\end{thm}
The following lemma shows that the classical nesting argument used
in the proof of Baire's theorem can be applied to $X$.

\begin{lemma}\label{Proposition 4.1.}
  Let $X$ be a locally compact space. Then
  for each $x\in X$ and each neighborhood $V\ni x$ there is a compact neighborhood $U\subset V$ of $x$
  ($X$ is then said to be microcompact).
\end{lemma}
\proof Let us denote by ${\frak V}_x$ the set of all neighborhoods
of $x$. Recall that a compact space is {\it regular}, i.e. for
each point $x\in X$ and each $V \in {\frak V}_x$, there is a
closed set $F \in {\frak V}_x$ with $F \subset V$. Let $V\in
{\frak V}_x$, and $K$ be a compact in ${\frak V}_x$. Then clearly
$V \cap K \in {\frak V}^{K}_x$, i.e. is a neighborhood of $x$ in
$K$. So, there is $F \in {\frak V}^{K}_x$ a closed set of $K$ with
$F \subset V \cap K$. So $F$ is compact, contained in $V$ and it
is easy to check that $F \in {\frak V}_x$.
\endproof
Using the classical nesting argument, this implies \ref{thm3}.

\medskip
{\it Remarks.} Locally compact cannot be weakened to locally
quasi-compact, since a countably infinite set with the finite
complement topology is a quasi-compact space which is not Baire.
This space is in fact microquasi-compact, that is, has the
property given in Lemma \ref{Proposition 4.1.} with `compact'
replaced by `quasi-compact'. The following chart summarizes the
relations between the local compactness properties, and how they
stand with respect to Baire (broken arrows mean ``does not
imply'').

\begin{figure}[h]
\centering
    \includegraphics[angle=0]{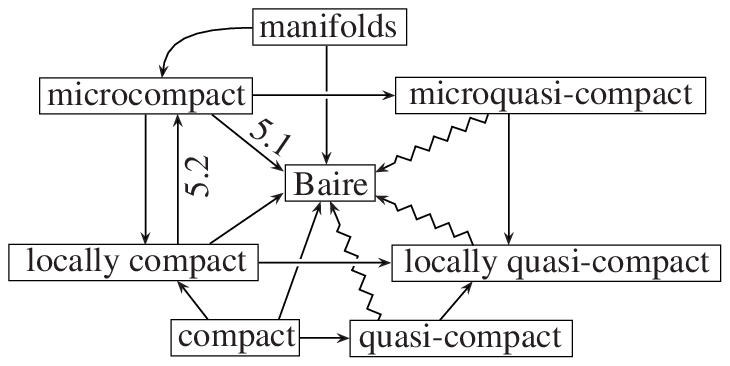}
\end{figure}

\bigskip \noindent{\bf Acknowledgments.} We thank A. Haefliger and
C. Weber for fruitful discussions.

\bigskip

\newbox\adress
\setbox\adress\vtop{\hsize 14.5 cm \noindent

\footnotesize
\noindent \textsc{Universit\'e de Gen\`eve, Section de
Math\'ematiques, 2-4, rue du Li\`evre, CP 240, 1211 Gen\`eve 24,
Suisse}

\noindent{\it E-mail address:} \verb"baillif@math.unige.ch"

\medskip
\noindent \textsc{Universit\'e de Gen\`eve, Section de
Math\'ematiques, 2-4, rue du Li\`evre, CP 240, 1211 Gen\`eve 24,
Suisse}

\noindent{\it E-mail address:} \verb"alexandregabard@hotmail.com"}

\vskip -0.0 cm

\centerline{\hbox{\copy\adress}}

\end{document}